\documentclass[a4paper,10pt]{academia}

\usepackage{graphicx}
\usepackage{amsfonts}
\usepackage{amssymb}
\usepackage{amsmath}
\usepackage{url}
\newtheorem{theorem}{Theorem}[section]

\newtheorem{definition}[theorem]{Definition}
\newtheorem{remark}[theorem]{Remark}

\def\RR{\mathbb R}

\def\CC{\mathbb C}
\def\pmatrix{ \left( \begin{array} }
\def\endpmatrix{ \end{array} \right) }

\def\no{\noindent}

\def\d{{\rm d}}

\def\pmatrix{ \left( \begin{array} }
\def\endpmatrix{ \end{array} \right) }

\def\cc{\gamma}
\def\dd{\delta}

\def\bfp{\boldsymbol{p}}

\def\bfp{\boldsymbol{p}}
\def\bfq{\boldsymbol{q}}

\def\bfy{\boldsymbol{y}}

\def\d2dxx{\frac{\partial^2}{\partial x^2}}

\def\no{\noindent}

\def\phi{\varphi}

\def\md{\mathrm{d}}

\title{Numerical comparisons between Gauss-Legendre methods and
Hamiltonian BVMs defined over Gauss points\thanks{Work developed within the
project ``Numerical methods and software for differential equations''.}}
\author{Luigi Brugnano$^\dagger$ \quad Felice Iavernaro$^\ddagger$ \quad Tiziana
Susca$^\ddagger$\\~\\
 $^\dagger$ \small Dipartimento di Matematica ``U.\,Dini'', Universit\`a di
Firenze, Italy\\
$^\ddagger$ \small Dipartimento di Matematica, Universit\`a di Bari, Italy}

\date{~\\ \em Dedicated to Prof.\,Manuel Calvo, on the occasion of his 65$th$ birthday.}

\begin{document}

\maketitle
\begin{abstract}
Hamiltonian Boundary Value Methods are a new class of energy
preserving one step methods for the solution of polynomial
Hamiltonian dynamical systems. They can be thought of as a
generalization of collocation methods in that they may be defined
by imposing a suitable set of {\em extended collocation
conditions}. In particular, in the way they are described in this
note, they are related to Gauss collocation methods with the
difference that they are able to precisely conserve the
Hamiltonian function in the case where this is a polynomial of any
high degree in the momenta and in the generalized coordinates. A
description of these new formulas is followed by a few test
problems showing  how, in many  relevant situations, the precise
conservation of the Hamiltonian is  crucial to simulate on a
computer the correct behavior of the  theoretical solutions.
\end{abstract}

\section{Introduction}
Hamiltonian Boundary Value Methods (HBVMs) form a subclass of
Boundary Value Methods (BVMs), whose main feature is that of
precisely conserving the Hamiltonian function associated with a
canonical  Hamiltonian system
\begin{equation}\label{hamilode}
\left\{  \begin{array}{l} \dot y =  J\nabla H(y),  \\
y(t_0) = y_0 \in\RR^{2m}, \end{array} \right.
 \qquad J=\pmatrix{rr} 0 & I \\ -I & 0 \endpmatrix \in \RR^{2m \times
 2m},
\end{equation}
($I$ is the identity matrix of dimension $m$),  in the case where
such function is of polynomial type.

Two key ideas have permitted the realization of HBVMs: the
definition of {\em discrete line integral} and what we called {\em
extended collocation conditions}. The former, first introduced in
\cite{IT0,IT1}, represents the discrete counterpart of the line
integral defined over conservative vector fields, while the second
is a relaxation of the classical collocation conditions which
assures the conservation of the energy along the numerical
solution $\{y_n\}$ generated by the method itself.

Just as an initial clarification, we briefly show how this new
approach to the problem reads when the classical Gauss collocation
method is considered (see \cite[Remark 2.1]{IT3} for more
details). Given a stepsize $h>0$ and a set of $s$ abscissae $c_1
<\dots <c_s$ disposed according to a Gauss-Legendre distribution
on $[0,1]$, the Gauss method of order $2s$ is defined by means of
the following polynomial collocation problem:
\begin{equation}
\label{collocation} \left\{ \begin{array}{ll}
\sigma(t_0)=y_0, \\[.2cm]
\dot \sigma(t_0+c_i h)= J \nabla H(\sigma(t_0+c_i h)), &i=1,\dots,s.
\end{array}
\right.
\end{equation}
As is well known, conditions \eqref{collocation} uniquely define a
polynomial $\sigma(t)$ of degree $s$ which is used to advance the
solution by posing $y_1= \sigma(t_0+h)$, while the internal stages
satisfy $Y_i=\sigma(t_0+c_i h)$, $i=1,\dots s$. The coefficients of
the Butcher array and the weights are given by
$$
b_j=\int_0^1 \ell_j(c) \mathrm{d} c, \qquad  a_{ij}=\int_0^{c_i}
\ell_j(c) \mathrm{d} c, \qquad \mbox{with }~ \ell_j(c) = \prod_{r
\not = j} \frac{c-c_r}{c_j-c_r}.
$$
The  $s$-degree  polynomial $\sigma(t)$ may be thought of as a
path in the phase space linking the state vectors $y_0$ to $y_1$
and passing through the stages $\{Y_i\}$. Due to the conservative
nature of the vector field, we have that
\begin{equation}
\label{1} H(y_1)-H(y_0) = \int_\sigma  \nabla H(y) \cdot \md y
=h\int_{0}^{1} \dot \sigma(t_0+\tau h)^T \nabla H(\sigma(t_0+\tau h))
\mathrm{d}\tau.
\end{equation}
Now, the above integral is exactly computed by the Gauss
quadrature formula with abscissae $\{c_i\}$ and weights $\{b_i\}$
if the degree of the integrand is not greater than $2s-1$ which
means that the degree of $H(y)$, say $\nu$, must not exceed $2$
(linear or quadratic Hamiltonians only). Under this assumption,
taking into account the collocation conditions
\eqref{collocation}, we obtain
\begin{equation}
\label{2} H(y_1)-H(y_0) =  h\sum_{i=1}^s b_i(\dot \sigma(t_i))^T
\nabla H(\gamma(t_i)) =  - h\sum_{i=1}^s b_i \nabla^T
H(\sigma(t_i)) J \nabla H(\gamma(t_i)) = 0,
\end{equation}
where $t_i=t_0+c_ih$. Thus,  by following a different route,  we
have obtained the classical result that the Gauss methods conserve
quadratic Hamiltonian functions while fails to conserve polynomial
Hamiltonian functions of higher degree.\footnote{\,This argument
may be generalized to other classes of collocation methods.}

The above example is the starting point of our approach: the {\em
discrete line integral} is the first sum in \eqref{2}, which turns
out to vanish for quadratic Hamiltonians, due to the collocation
conditions \eqref{collocation}.


The next section reports a descriptive introduction to HBVMs with
much emphasis to the key ideas they rely on. We refer the reader to
the papers \cite{BIT, BIT1, IT3, BIT2, BIT3, IP1, IP2, BIS} for the
details about the basic theory and implementation of HBVMs, and to
the monograph \cite{BT1} as a reference for the theory of BVMs.

In Section \ref{test_sec} we report a number of test problems of
some relevance in the literature, for which the precise
conservation of the energy turns out to be a crucial feature for
the  correct reproduction of the long time behavior of the
solutions. This will be testified by comparing  HBVMs to Gauss
methods which, by the way, are symplectic integrators.

\section{Hamiltonian Boundary Value Methods}\label{hbvm_sec}
In this section we introduce HBVMs by slightly elaborating the
arguments in \cite{BIT,BIT1,BIT3}. As said above, the basic idea
which HBVMs rely on is the so called discrete line integral, which
is the discrete counterpart of the line integral associated with a
conservative vector field. In more detail, starting from
(\ref{1}), we consider a polynomial, of degree at most $s$, such
that
\begin{equation}\label{01}\sigma(t_0) = y_0, \qquad \sigma(t_0+h)
= y_1,\end{equation}

\no providing an approximation to the solution on the interval
$[t_0,t_0+h]$. We consider the following expansions,
\begin{equation}\label{s1s} \dot\sigma(t_0+ch) = \sum_{j=1}^s P_j(c)\cc_j, \qquad
\sigma(t_0+ch) = y_0 + h\sum_{j=1}^s \cc_j \int_{0}^c P_j(\tau)\md
\tau,\end{equation}

\no where the (vector) coefficients $\{\cc_i\}$ are to be
determined. We also assume that the polynomials $\{P_i\}$
constitute an orthonormal basis, on the interval $[0,1]$, for the
vector space $\Pi_{s-1}$ of polynomials of degree at most $s-1$,
i.e.,
$$ \int_0^1 P_i(\tau)P_j(\tau) = \dd_{ij}, \qquad
i,j=1,\dots,s,$$

\no with $\dd_{ij}$ the Kronecker symbol. Such polynomials can be
easily obtained by a suitable scaling of the shifted Legendre
polynomials \cite{BIT3}. Substituting the first expansion in
(\ref{s1s}) into the line integral in (\ref{1}), and requiring the
resulting expression to vanish, then gives

$$\sum_{j=1}^s \cc_j^T\int_0^1 P_j(\tau)\nabla H(\sigma(t_0+\tau h))\md\tau=0,$$

\no which is certainly satisfied by choosing

\begin{equation}\label{ccj}
\cc_j = \int_0^1 P_j(\tau)J\nabla H(\sigma(t_0+\tau h))\md\tau, \qquad
j=1,\dots,s.
\end{equation}

\no Multiplication of (\ref{ccj}) by $h\int_0^c P_j(x)\md x$ and
summation over $j$ then gives, by virtue of the second expansion
in (\ref{s1s}),

\begin{equation}\label{infhbvm}
\sigma(t_0+c h) = y_0 + h\sum_{j=1}^s \int_0^c P_j(x)\md x \,
\int_0^1 P_j(\tau)J\nabla H(\sigma(t_0+\tau h))\md\tau, \qquad
c\in[0,1].\end{equation}

Let us now assume that $H(y)$ is a polynomial of degree $\nu$.
Consequently, the integral appearing at the right-hand side in
(\ref{infhbvm}) can be exactly discretized by a Gaussian formula
over $k$ Gauss-Legendre abscissae $\{c_i\}$, which we shall consider
hereafter, provided that
\begin{equation}\label{kk} k\ge\frac{\nu s}2.\end{equation}

\no Let us denote by $\{\omega_i\}$ the weights of the quadrature
formula in the interval $[0,1]$, and
set\begin{equation}\label{yiaij}y_i = \sigma(t_0+c_ih), \quad
a_{ij} = \int_0^{c_i}P_j(x)\md x, \qquad i=1,\dots,k, \quad
j=1,\dots,s.\end{equation}

\no Consequently, (\ref{infhbvm}) can be (exactly) discretized as:
\begin{equation}\label{hbvmks}
y_i = y_0 + h\sum_{j=1}^s a_{ij} \sum_{\ell=1}^k \omega_\ell
P_j(c_\ell)J\nabla H(y_\ell), \qquad i=1,\dots,k.
\end{equation}

\begin{definition}\label{hbvmdef} The set of equations
(\ref{hbvmks}), to be solved for the unknowns $\{y_i\}$, defines
an {\em HBVM with $k$ steps and degree $s$}, in short
HBVM$(k,s)$.\end{definition}

\no For such a method, the following properties hold true
\cite{BIT1}:
\begin{itemize}

\item it has order $2s$ for all $k\ge s$;

\item it is symmetric and perfectly $A$-stable (i.e., its
stability region coincides with the left-half complex plane,
$\CC^-$ \cite{BT1});

\item for $k=s$, it reduces to the Gauss-Legendre method of order $2s$;

\item it exactly preserves polynomial Hamiltonian functions of degree $\nu$,
provided that (\ref{kk}) holds true.

\end{itemize}

\begin{remark} The actual implementation of HBVM$(k,s)$ can be seen to
result in the solution of a system of (block) size $s$, whatever
is the value of $k$ considered \cite{BIT,BIT3}. Consequently, if
needed, large values of $k$ can be easily considered.\end{remark}

The arguments in the previous remark, allow us to consider the limit
formula of (\ref{yiaij})--(\ref{hbvmks}), in the case where $H(y)$
is non-polynomial, as $k\rightarrow\infty$. Clearly such a limit is
given by formula (\ref{infhbvm}), which, according to \cite{BIT1},
is named {\em HBVM$(\infty,s)$} or {\em $\infty$-HBVM of degree
$s$}.

However, we emphasize that formula (\ref{infhbvm}) becomes an
operative method only after that a suitable discretization of the
inner integral is considered and, replacing the integral by a
quadrature formula with $k$ nodes, leads back to a HBVM$(k,s)$
method.

One can easily argue that, since in the non polynomial case the
quadrature formula can approximate  the corresponding integral
with an arbitrary accuracy, under suitable regularity assumptions
for $H(y)$, a {\em practical} conservation of the energy may be
obtained \cite{BIT1,IT2}. The term ``practical'' means that, in
many general situations, when $k$ is high enough, the method makes
no distinction between the function $H(y)$ and its polynomial
approximation, being the latter in a neighborhood of size
$\varepsilon$ of the former, where $\varepsilon$ denotes the
machine precision.

We end this section by observing that, by differentiating both
members of (\ref{infhbvm}), one obtains
$$\dot\sigma(t_0+c h) = \sum_{j=1}^s P_j(c) \,
\int_0^1 P_j(\tau)J\nabla H(\sigma(t_0+\tau h))\md\tau, \qquad
c\in[0,1],$$

\no which at the points $\{c_i\}$ provides, assuming $H(y)$ to be
a polynomial and $k$ large enough:
$$\dot\sigma(t_0+c_i h) = \sum_{j=1}^s P_j(c_i) \,
\int_0^1 P_j(\tau)J\nabla H(\sigma(t_0+\tau h))\md\tau, \qquad
i=1,\dots,k.$$

\no Such formulae (the former being the limit of the latter as
$k\rightarrow\infty$) can be regarded as a kind of {\em extended
collocation conditions} that generalize conditions
\eqref{collocation}, according to \cite[Section\,2]{IT3} (see also
\cite{BIT1}).

\section{Numerical tests}
\label{test_sec} We present a few numerical test highlighting the
good behavior of HBVMs in the long-time simulation of Hamiltonian
systems. A direct comparison of HBVMs with Gauss methods  is
reported in order to better emphasize the stability properties of
the former methods even when compared to a well known class of
symplectic formulae.\footnote{As was seen in the previous section,
the choice of Gauss methods has also been dictated by the fact
that they represent the {\em generating} formulae of HBVMs when we
use a Gauss distribution of the abscissae, namely the Gauss method
of order $2s$ coincides with HBVM($s$,$s$).}

The use of a large stepsize of integration is a prerogative in
long-time simulation of an evolutionary problem but, in general,
one is forced to reduce $h$ under a critical threshold in order to
guarantee the qualitative behavior of the theoretical solution to
be well reproduced by the numerical solution. From this point of
view, we show that HBVMs allow the use of larger stepsizes than
Gauss methods, which states that the conservation of the
Hamiltonian function plays an important role in detecting the
correct topological features of the solutions.

\subsection{Sitnikov's problem}
One of the main problems in Celestial Mechanics is to describe the
motion of $N$ point particles of positive mass $\{m_i\}$ moving
under Newton's law of gravitation when we know their positions
$\{q_i\}$ and momenta $\{p_i\}$ at a given time. Such a dynamical
system, called the $N$-body problem, is in the form
\eqref{hamilode}, with Hamiltonian
\begin{equation}
\label{kepler} H(\bfq,\bfp)=\frac{1}{2} \sum_{i=1}^N
\frac{||p_i||_2^2}{m_i} - G \sum_{i=1}^N m_i
\sum_{j=1}^{i-1}\frac{m_j}{||q_i-q_j||_2},
\end{equation}
with $G$ the gravitational constant. While the two-body problem is
completely solved in the sense that we can describe explicitly all
its solutions (see, e.g., \cite{HLW}), this is no more the case,
for $N\ge 3$. Consequently, numerical simulation is of interest,
in such a case.

The Sitnikov problem is a particular configuration of the $3$-body
dynamics. In this problem two bodies of equal mass (primaries)
revolve about their center of mass, here placed at the origin, in
elliptic orbits in the $(x,y)$-plane. A third, and much smaller
body (planetoid), is placed on the $z$-axis with initial velocity
parallel to this axis as well.

\begin{figure}[ht]
\begin{center}
\includegraphics[width=7.7cm,height=6cm]{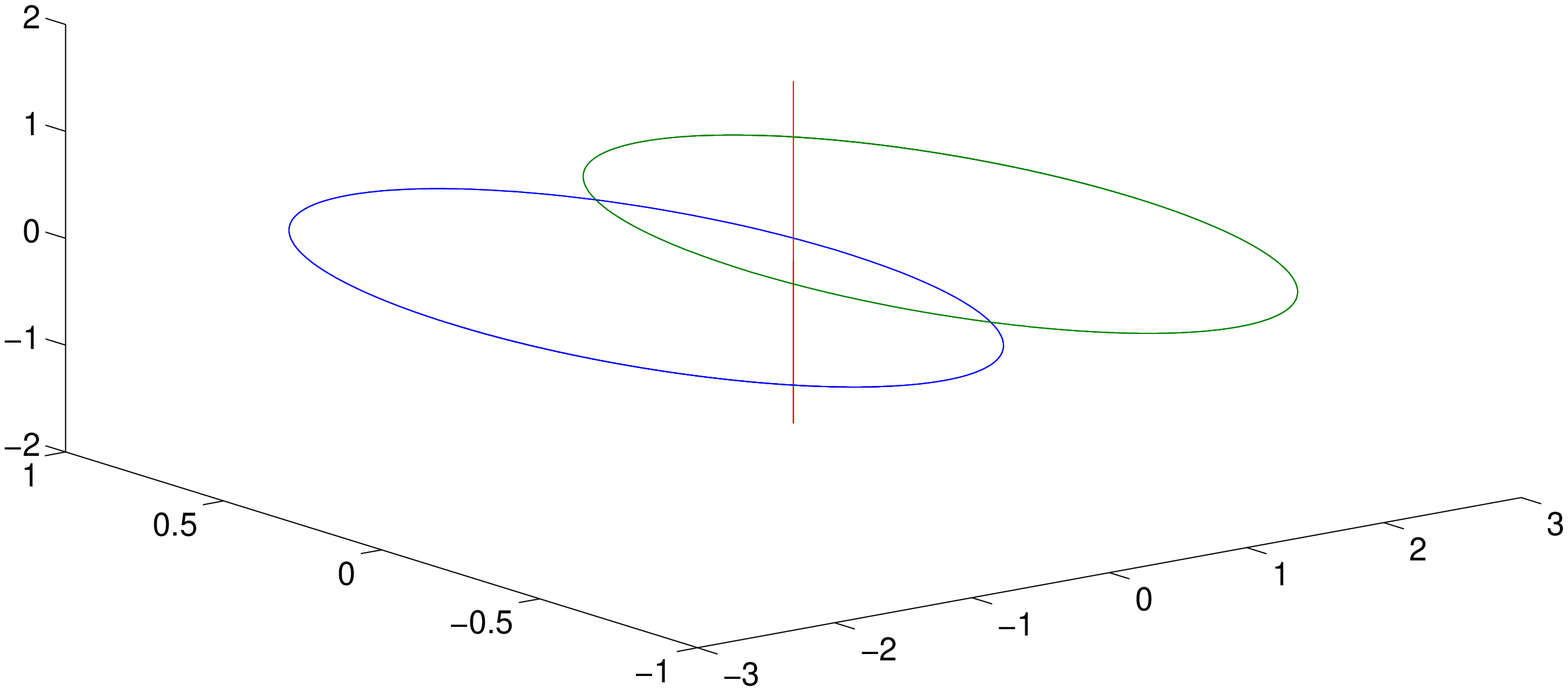}\hspace{.2cm}
\includegraphics[width=7.7cm,height=6cm]{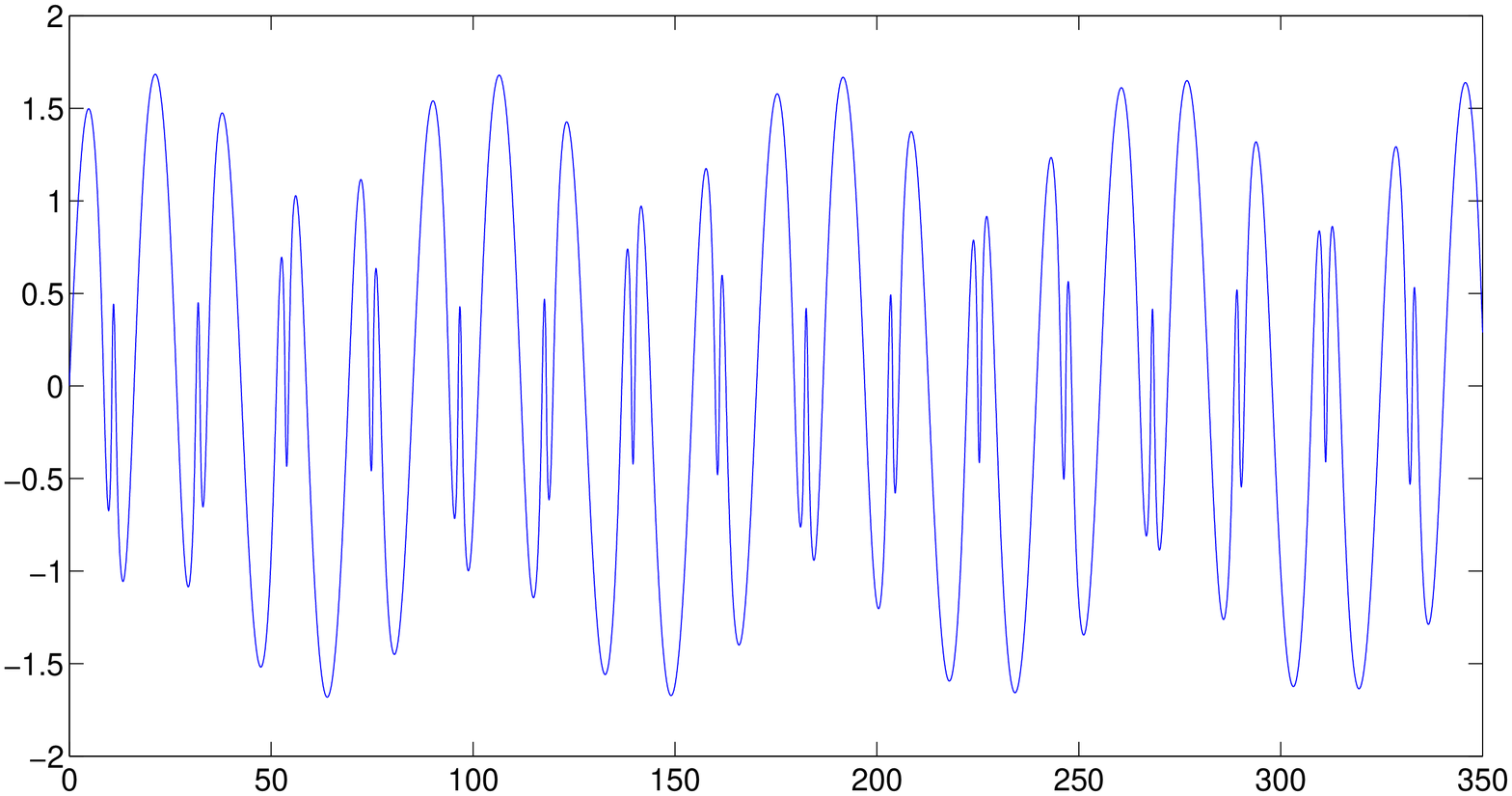}
\end{center}
\vspace*{-0.7cm}  \caption{The left picture displays the
configuration of  $3$-bodies in the Sitnikov problem. To an
eccentricity of the orbits of the primaries $e=0.75$, there
correspond bounded chaotic oscillations of the planetoid as is
argued by looking at the space-time diagram in the right picture.}
\label{sit_fig1}
\begin{center}
\includegraphics[width=7.7cm,height=6cm]{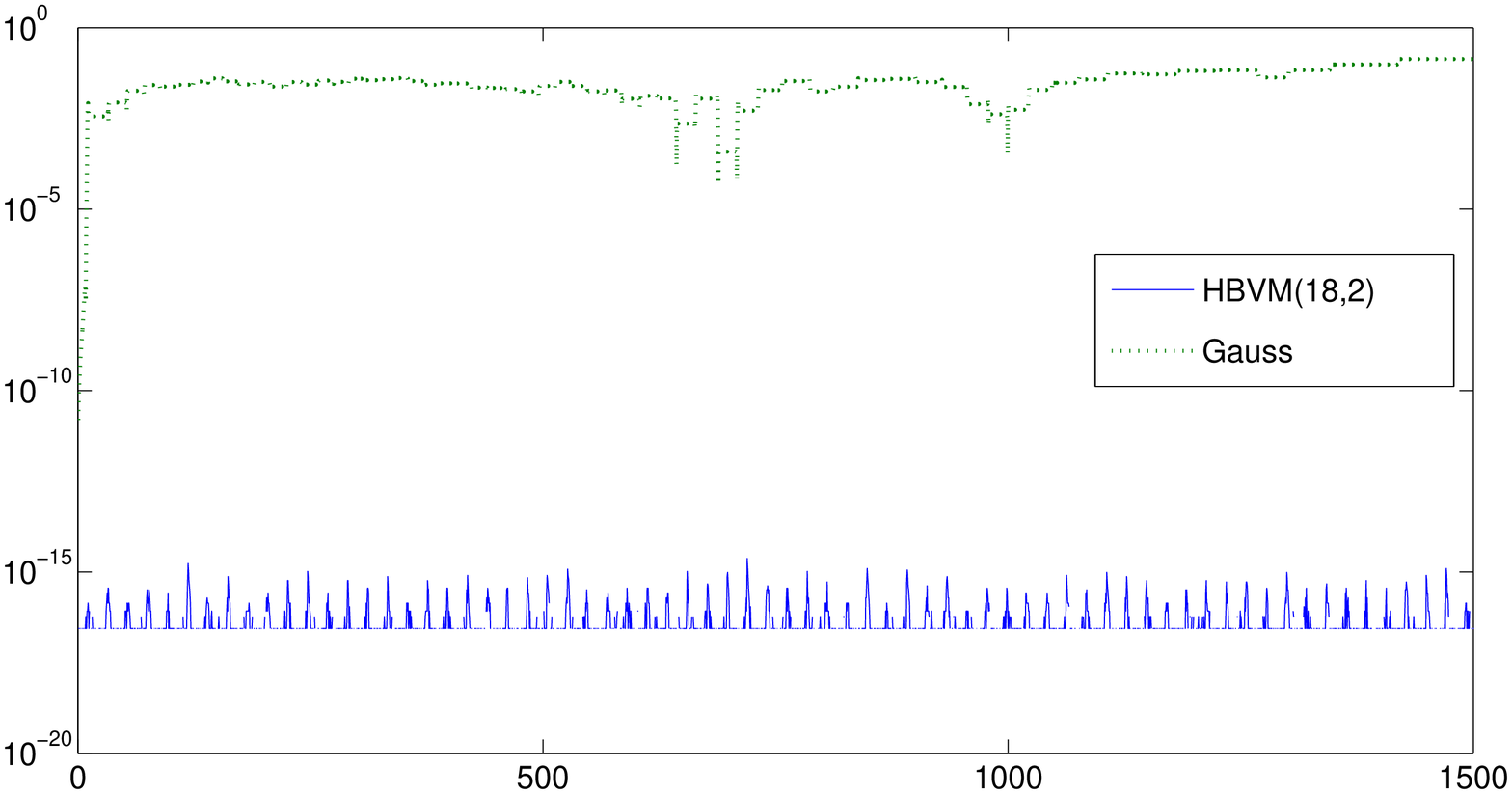}
\hspace*{0.0cm}
\includegraphics[width=7.7cm,height=6cm]{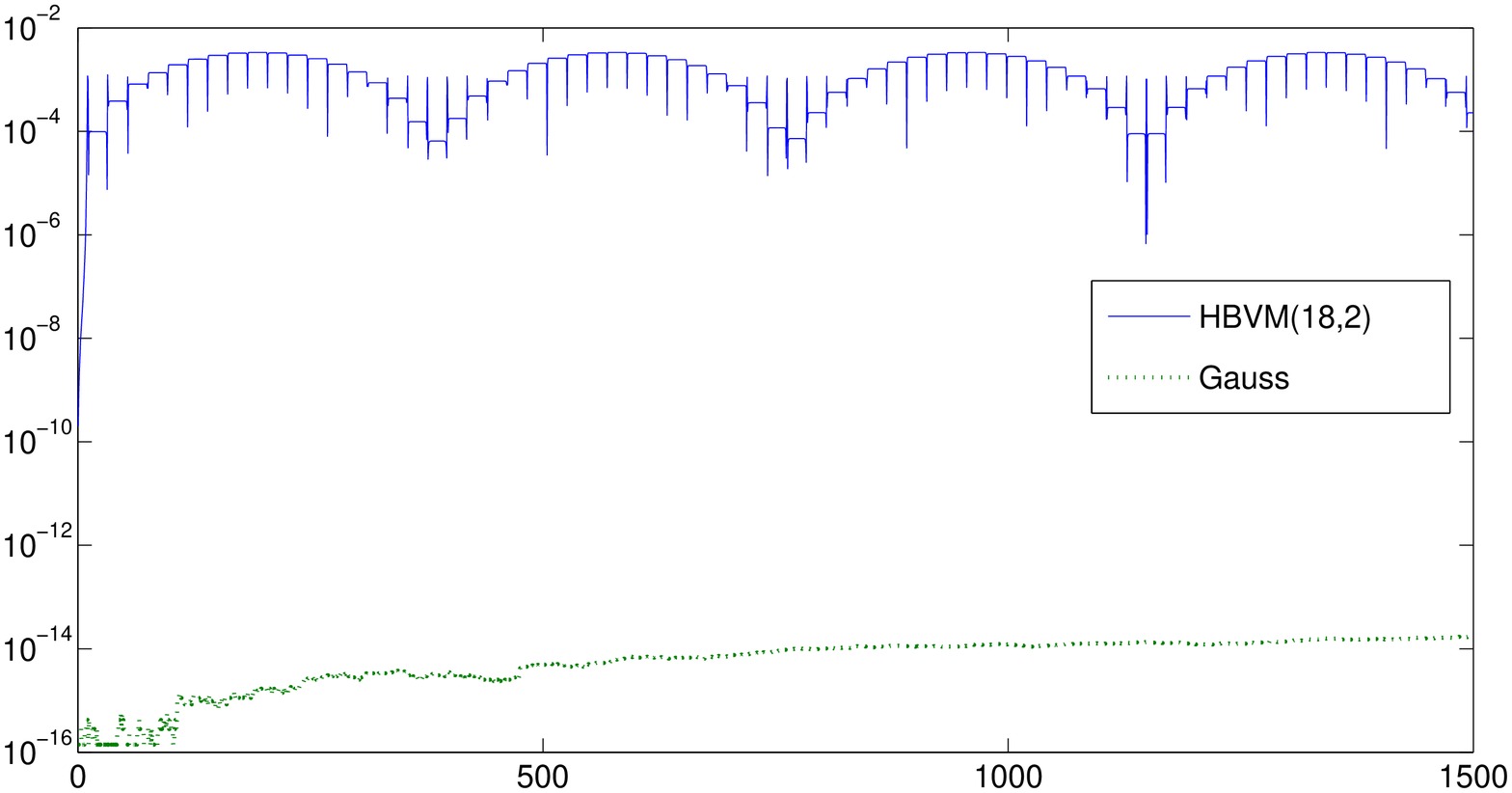}
\end{center}
\vspace*{-0.7cm} \caption{Left picture: relative error
$|H(y_n)-H(y_0)|/|H(y_0)|$ of the Hamiltonian function evaluated
along the numerical solution of the HBVM($18$,$2$) and the Gauss
method. Right picture: relative error $|M(y_n)-M(y_0)|/|M(y_0)|$
of the angular momentum evaluated along the numerical solution of
the HBVM($18$,$2$) and the Gauss method.} \label{sit_fig5}
\end{figure}

The third body is small enough that the two body dynamics of the
primaries is not destroyed. Then, the motion of the third body
will be restricted to the $z$-axis and oscillating around the
origin but not necessarily periodic. In fact, this problem has
been shown to exhibit a chaotic behavior  when the eccentricity of
the orbits of the primaries exceeds a critical value that, for the
data set we have used, is $\bar e \simeq 0.725$ (see Figure
\ref{sit_fig1}).

We have solved the Kepler problem with Hamiltonian function
\eqref{kepler} by the Gauss method of order 4 (HBVM(2,2)) and by
HBVM(18,2) (order 4 and $18$ steps), with the following set of
parameters:

\begin{center}
\begin{tabular}{ccccccccc}
$N$ & $G$ & $m_1$ & $m_2$ & $m_3$ & $e$ & $d$ & $h$ & $t_{\mbox{max}}$
\\[.1cm]
\hline \\[-.4cm] $3$ & $1$ & $1$ & $1$ & $10^{-5}$ & $0.75$ & $5$ & $0.5$ &
$1500$
\end{tabular}
\end{center}
where $e$ is the eccentricity, $d$ is the distance of the
apocentres of the primaries (points at which the two bodies are
the furthest), $h$ is the time-step and $[0,\,t_{\mbox{max}}]$ is
the time integration interval. The eccentricity $e$ and the
distance $d$ may be used to define the initial condition
$\bfy_0=[\bfq_0,\bfp_0]$ (see \cite{James} for the details):
$$
\begin{array}{l}
\bfq_0 = [-\frac{5}{2},~  0,~  0,~  \frac{5}{2},~  0,~  0,~   0,~
0,~10^{-9}],\\[.1cm]
\bfp_0 =[0,~ -\frac{1}{20}\sqrt{10},~  0,~    0,~
\frac{1}{20}\sqrt{10},~ 0,~ 0,~ 0,~ \frac{1}{2}].
\end{array}
$$

\begin{figure}[ht]
\begin{center}
\includegraphics[width=7.7cm,height=6cm]{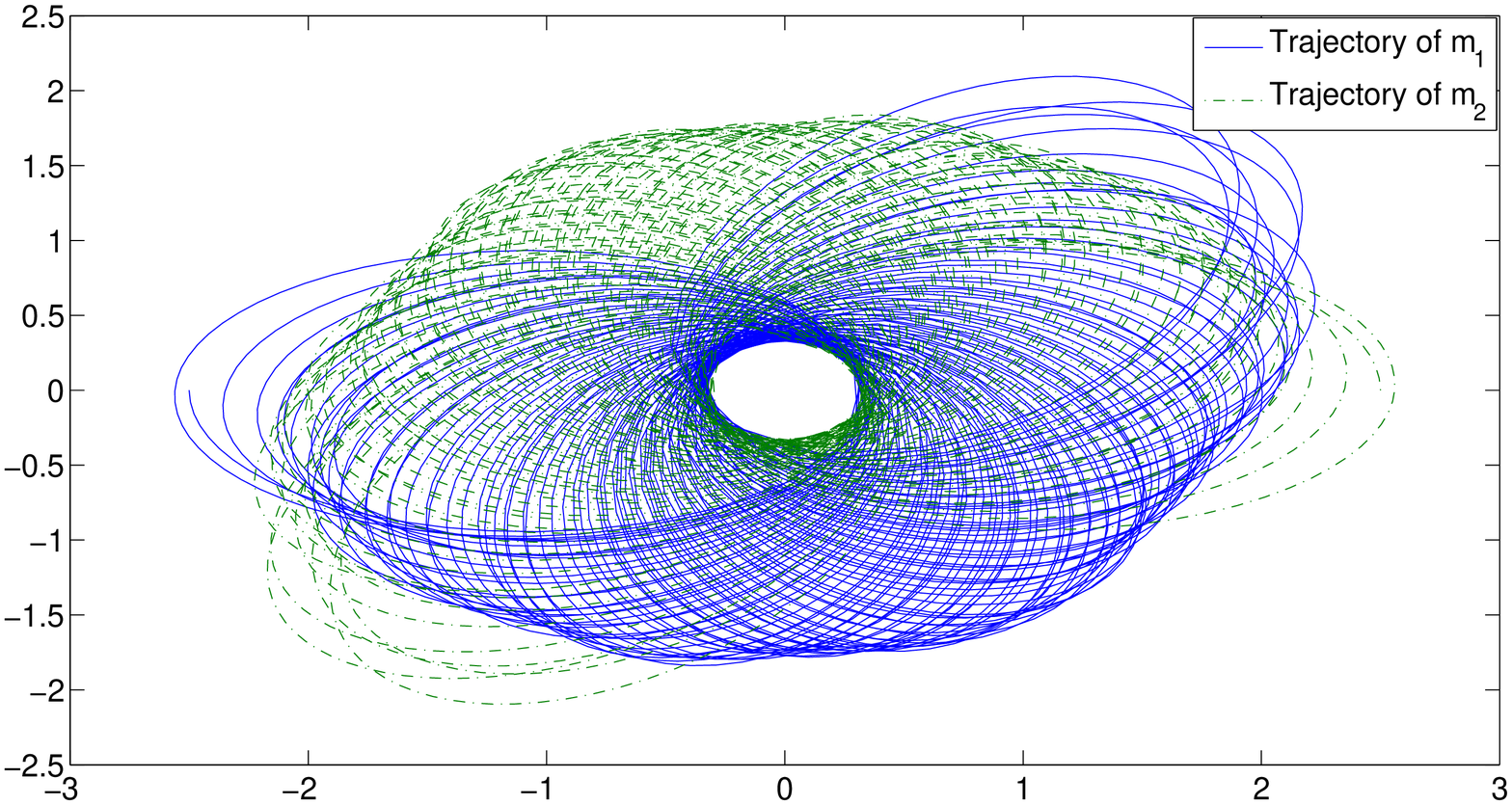}
\includegraphics[width=7.7cm,height=6cm]{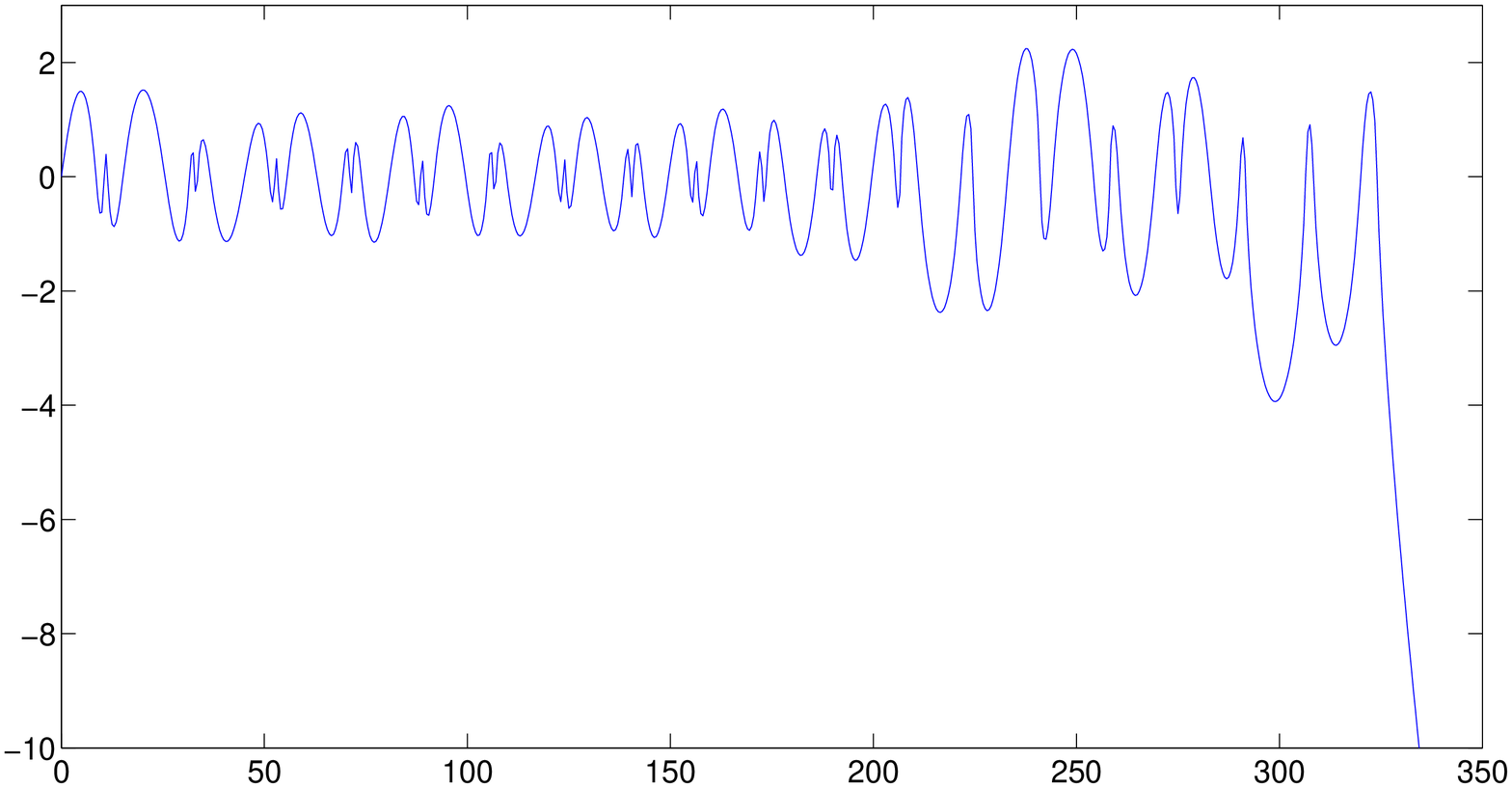}
\end{center}
\vspace*{-0.7cm}  \caption{The Sitnikov problem solved by the
Gauss method of order 4, with stepsize $h=0.5$, in the time
interval $[0,1500]$. The trajectories of the primaries in the
$(x,y)$-plane (left picture) exhibit a very irregular behavior
which causes the planetoid to eventually leave the system, as
illustrated by the space-time diagram in the right picture.}
\label{sit_fig2}
\begin{center}
\includegraphics[width=7.7cm,height=6cm]{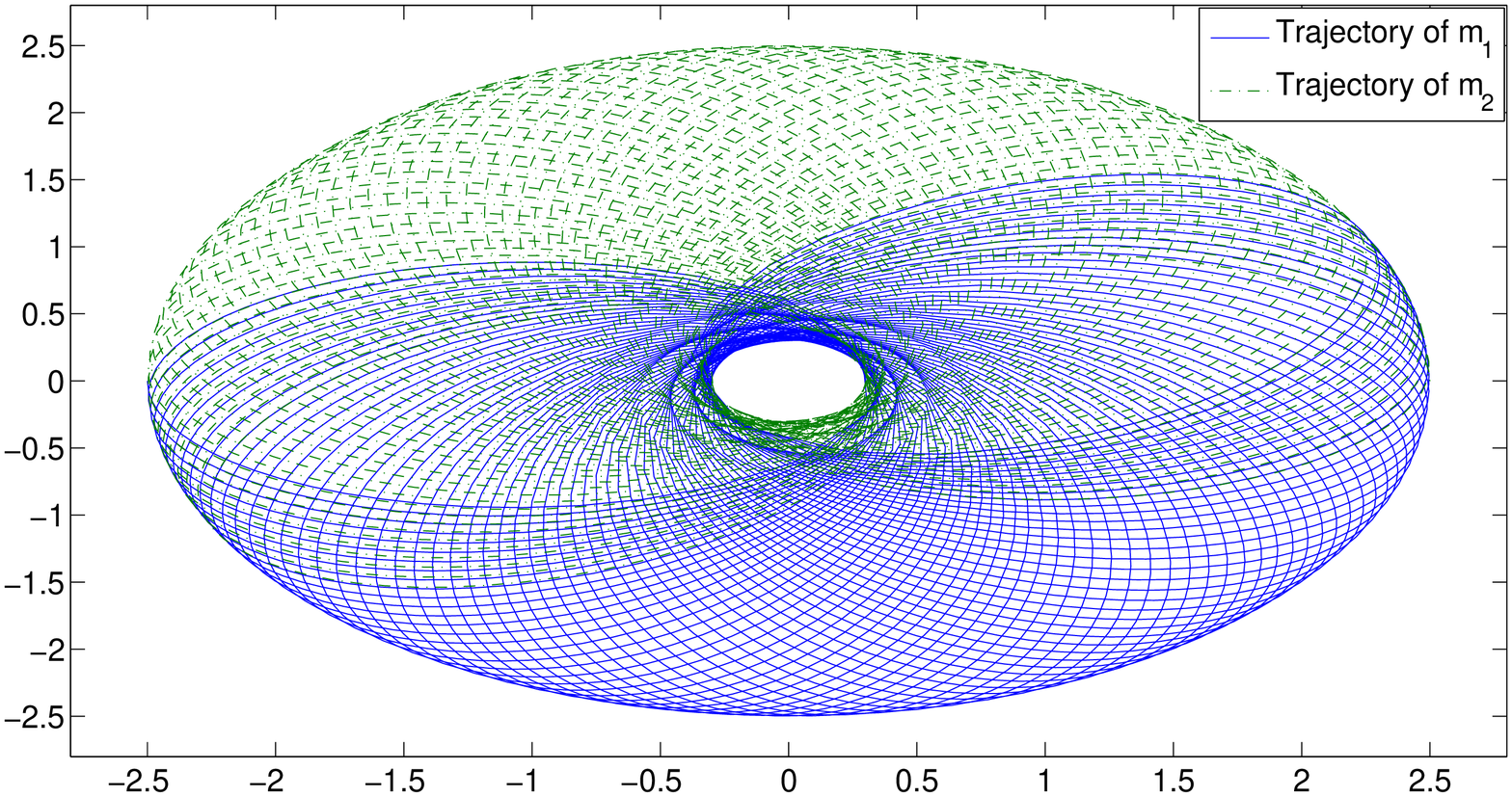}
\hspace*{0.0cm}
\includegraphics[width=7.7cm,height=6cm]{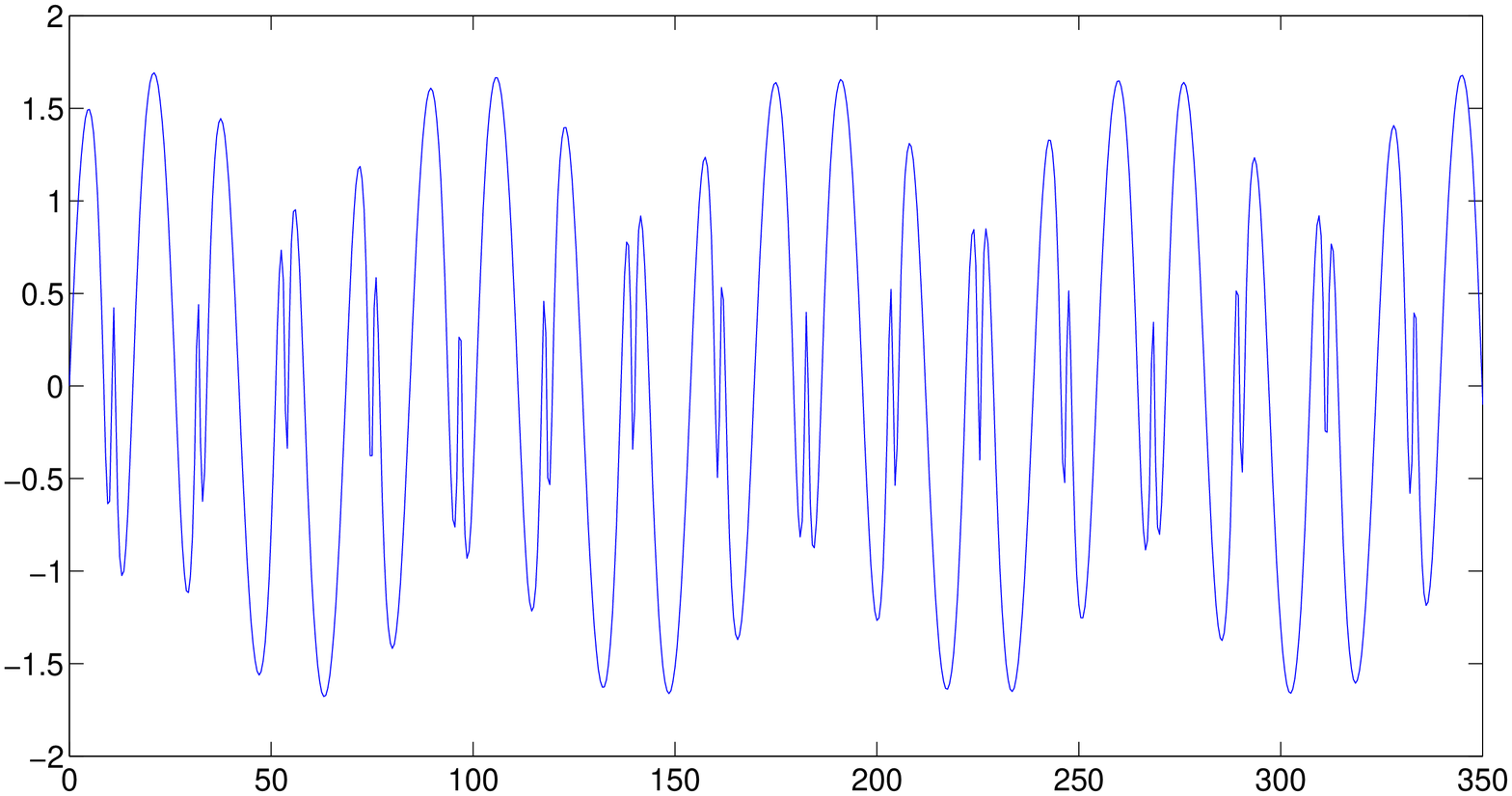}
\end{center}
\vspace*{-0.7cm} \caption{The Sitnikov problem solved by the
HBVM(18,2) method (order 4), with stepsize $h=0.5$, in the time
interval $[0,1500]$. Left picture: the trajectories of the
primaries are ellipse shape. The discretization introduces a
fictitious uniform rotation of the $(x,y)$-plane which, however,
does not alter the global symmetry of the system. Right picture:
the space-time diagram of the planetoid on the $z$-axis displayed
(for clearness) on the time interval $[0, 350]$ shows that,
although a large value of the stepsize $h$ has been used, the
overall behavior of the dynamics is well reproduced (compare with
the right picture of Figure \ref{sit_fig1}).} \label{sit_fig3}
\end{figure}

\begin{figure}[ht]
\begin{center}
\includegraphics[width=7.7cm,height=6cm]{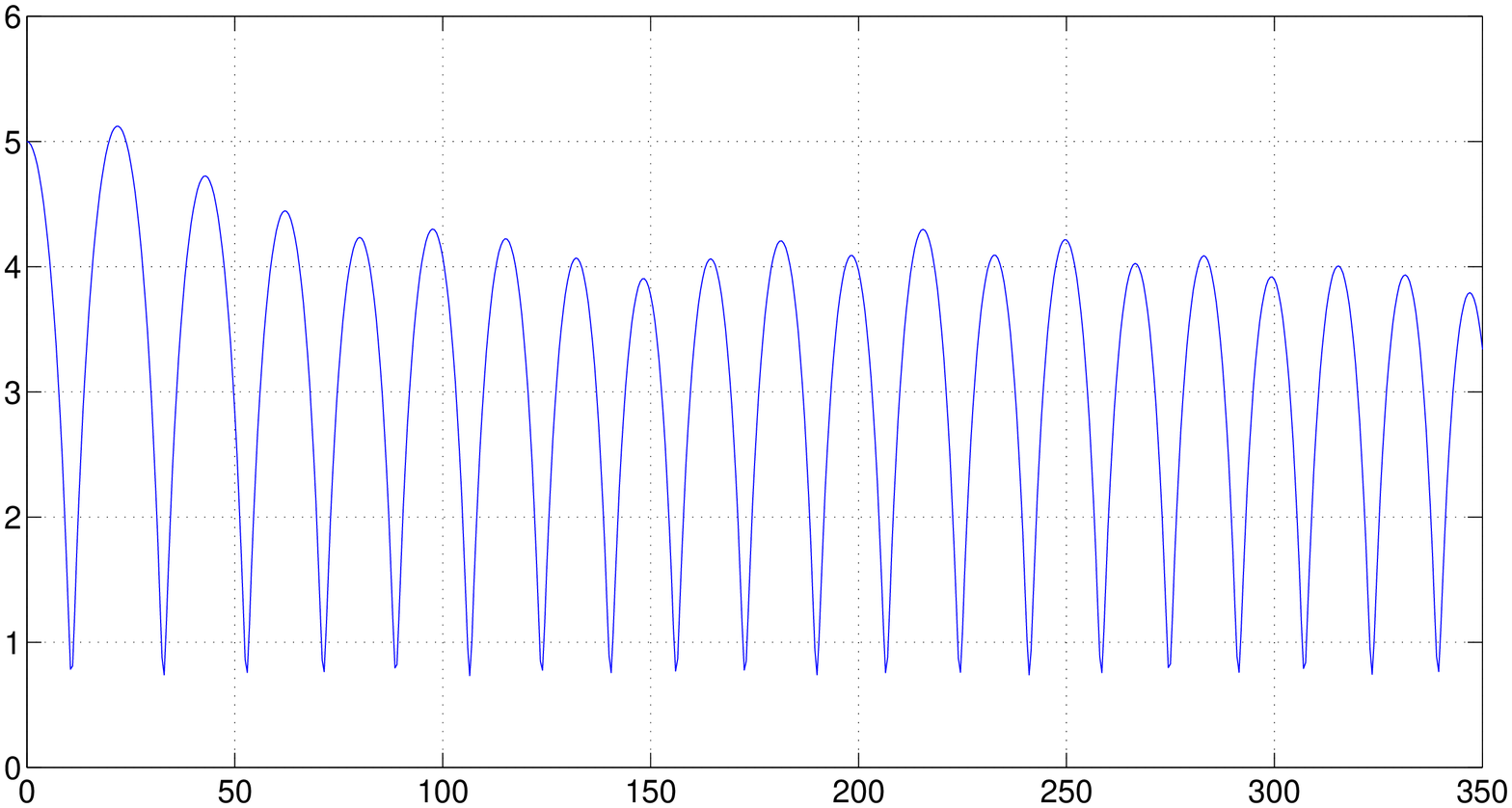}
\includegraphics[width=7.7cm,height=6cm]{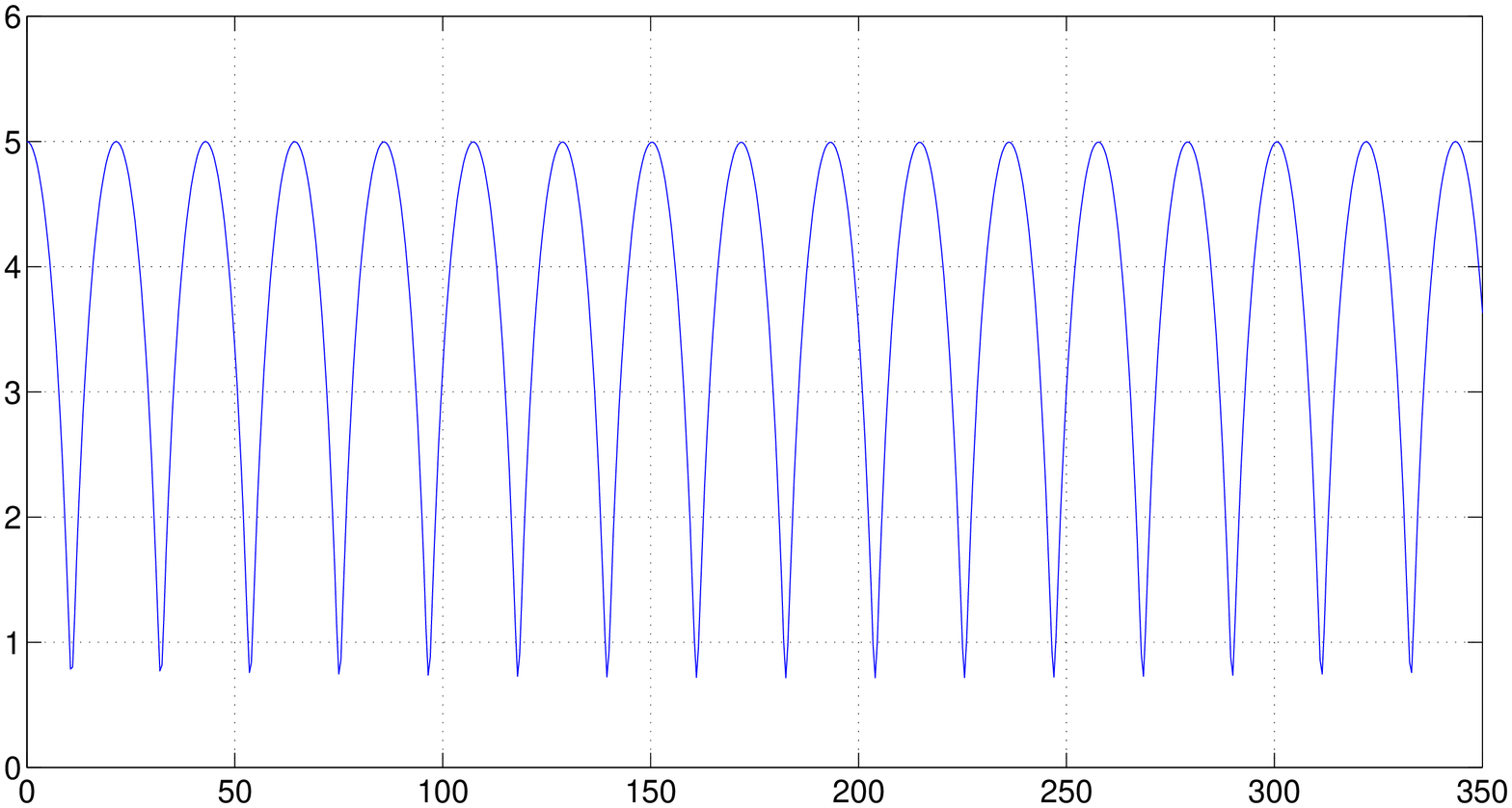}
\end{center}
\vspace*{-0.7cm}  \caption{Distance between the two primaries as a
function of the time, related to the numerical solutions generated
by the Gauss method (left picture) and HBVM(18,2) (right picture).
The maxima correspond to the distance of apocentres. These are
conserved by HBVM(18,2) while the Gauss method introduces patchy
oscillations that destroy the overall symmetry of the system.}
\label{sit_fig4}
\end{figure}

First of all, we consider the two pictures in Figure
\ref{sit_fig5} reporting the relative errors in the  Hamiltonian
function and in the angular momentum evaluated along the numerical
solutions computed by the two methods. According to (\ref{kk}), we
know that the HBVM(18,2) precisely conserves Hamiltonian
polynomial functions of degree at most $18$. This accuracy is high
enough to guarantee that the nonlinear Hamiltonian function
\eqref{kepler} is as well conserved up to the machine precision
(see the left picture): from a geometrical point of view, this
means that a local approximation of the level curves of
\eqref{kepler} by a polynomial of degree $18$ leads to a
negligible error. The Gauss method exhibits a certain error in the
Hamiltonian function while, being this formula symplectic, it
precisely conserves the angular momentum, as is confirmed by
looking at the right picture of Figure \ref{sit_fig5}. From the
same picture, one sees that the error in the numerical angular
momentum associated with the HBVM(18,2) undergoes some bounded
periodic-like oscillations.

Figures \ref{sit_fig2} and \ref{sit_fig3} show the numerical
solution computed by the Gauss method and HBVM(18,2),
respectively. Since the methods leave the $(x,y)$-plane invariant
for the motion of the primaries and the $z$-axis invariant for the
motion of the planetoid, we have just reported the motion of the
primaries in the $(x,y)$-phase plane (left pictures) and the
space-time diagram of the planetoid (right picture).

We observe that, for the Gauss method, the orbits of the primaries
are irregular in character so that the third body, after
performing some oscillations around the origin, will eventually
leave the system (see the right picture of Figure \ref{sit_fig2}).
On the contrary (left picture of Figure \ref{sit_fig3}), the
HBVM(18,2) generates a quite regular phase portrait. Due to the
large stepsize $h$ used, a sham rotation of the $(x,y)$-plane
appears which, however, does not destroy the global symmetry of
the dynamics, as testified by the bounded oscillations of the
planetoid (right picture of Figure \ref{sit_fig3}) which look very
similar to the reference ones in Figure \ref{sit_fig1}. This
aspect is also confirmed by the pictures in Figure \ref{sit_fig4},
displaying the distance of the primaries as a function of the
time. We see that the distance of the apocentres (corresponding to
the maxima in the plots), as the two bodies wheel around the
origin, are preserved by the HBVM(18,2) (right picture) while the
same is not true for the Gauss method (left picture).

\subsection{The H\'enon-Heiles problem}
The H\'{e}non-Heiles equation originates from a problem in
Celestial Mechanics describing the motion of a star under the
action of a gravitational potential of a galaxy which is assumed
time-independent and with an axis of symmetry (the $z$-axis) (see
\cite{HH} and references therein). The main question related to
this model was to state the existence of a third first integral,
beside the total energy and the angular momentum.\footnote{An
analytical approach to the problem may be found in \cite{Gu},
where the author finds out a formal expansion of the third
invariant.} By exploiting the symmetry of the system and the
conservation of the angular momentum, H\'enon and Heiles reduced
from three (cylindrical coordinates) to two (planar coordinates)
the degrees of freedom, thus showing that the problem was
equivalent to the study of the motion of a particle in a plane
subject to an arbitrary potential $U(q_1,q_2)$:
\begin{equation}
\label{HH} H(\bfq,\bfp)=\frac{1}{2}(p_{1}^2+p_{2}^2)+U(q_1,q_2).
\end{equation}

\begin{figure}[t]
\begin{center}
\includegraphics[width=13cm,height=8cm]{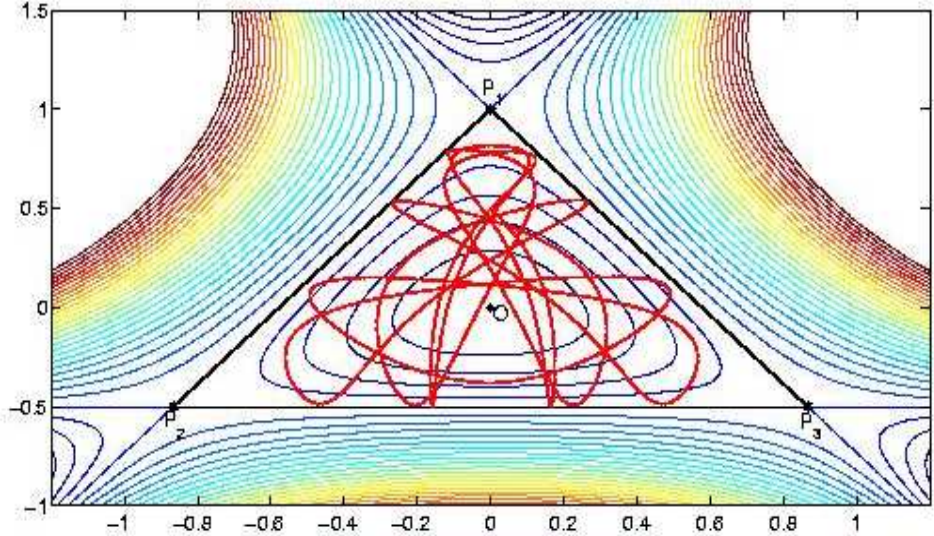}
\end{center}
\vspace*{-0.7cm} \caption{Level curves of the potential
$U(q_1,q_2)$ of the H\'enon-Heiles problem (see
\eqref{Henon_potential}). The origin $O$ is a stable equilibrium
point, whose domain of stability contains the equilateral triangle
having as vertices the saddle points $P_1$, $P_2$, and $P_3$,
provided that the total energy does not exceed the value
$\frac{1}{6}$. Inside the triangle an orbit $(q_1(t),q_2(t))$ is
traced whose total energy is close (but lower than) $\frac{1}{6}$.
The trajectory gets very close to the sides of the triangle, which
makes the problem of conserving the total energy in the numerical
solution an important feature to avoid instability when a large
stepsize is used.}\label{henon_fig1}
\end{figure}

Since $U$ in \eqref{HH} has no symmetry in general, we cannot
consider the angular momentum as an invariant anymore, so that the
only known first integral is the total energy represented by
\eqref{HH} itself, and the question is whether or not a second
integral does exist. H\'enon and Heiles conducted a series of tests
with the aim of giving a numerical evidence of the existence of such
integral for moderate values of the energy $H$, and of the
appearance of chaotic behavior when $H(\bfq,\bfp)$ becomes larger
than a critical value. In particular, for their experiments they
choose
\begin{equation}
\label{Henon_potential}
U(q_1,q_2)=\frac{1}{2}(q_{1}^2+q_{2}^2)+q_{1}^2q_{2}-\frac{1}{3}q_{2}^3,
\end{equation}
which makes the Hamiltonian function a polynomial of degree three.

When $U(q_1,q_2)$ approaches the value $\frac{1}{6}$, the level
curves of $U$ tend to an equilateral triangle, whose vertices are
saddle points of $U$ (see Figure \ref{henon_fig1}). This vertices
have coordinates $P_1=(0,1)$,
$P_2=(-\frac{\sqrt{3}}{2},-\frac{1}{2})$ and
$P_3=(\frac{\sqrt{3}}{2},-\frac{1}{2})$.

We consider an initial point $(\bfq_0,\bfp_0)$ such that
$\bfq_{0}$ is inside the triangle $U\leq\frac{1}{6}$ and
$H(\bfq_{0},\bfp_{0})<\frac{1}{6}$: then the orbit originating
from $(\bfq_0,\bfp_0)$ will never abandon the triangle for any
value of the time $t$. However, when $H(\bfq_{0},\bfp_{0})$ is
chosen very close to $\frac{1}{6}$, a numerical method which does
not preserve exactly the total energy could cause the (numerical)
orbit to jump outside the triangle and possibly to diverge to
infinity. This aspect is further emphasized when a large stepsize
of integration is used, as is usually required in the long time
simulation of a dynamical system.

\begin{figure}[ht]
\begin{center}
\includegraphics[width=13cm,height=8cm]{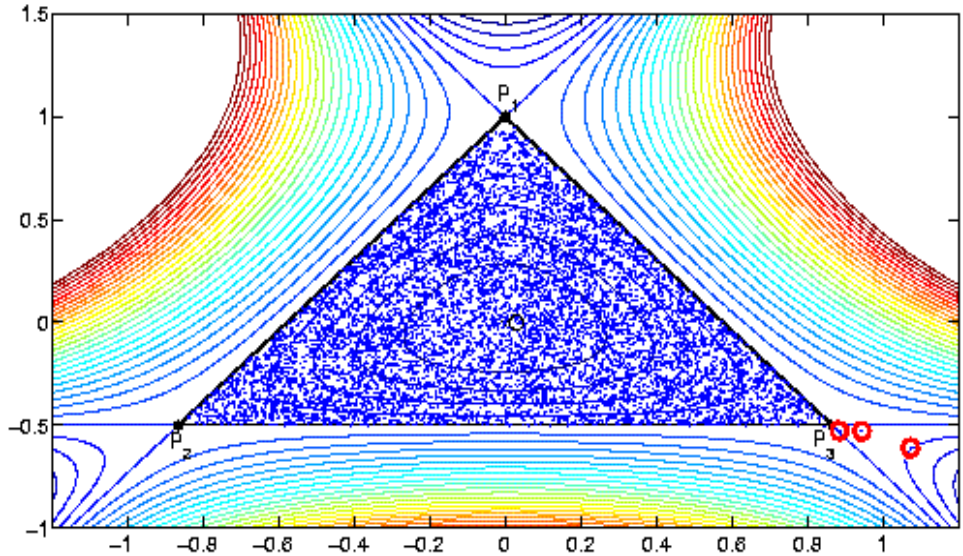}
\end{center}
\vspace*{-0.7cm}  \caption{The numerical trajectory in the
$(q_1,q_2)$-plane computed by the Gauss method of order four with
stepsize $h=1$. The stable character of the continuous orbit is
not correctly reproduced by the numerical method: after a time
$t\simeq 7000$ the orbit escapes from the triangle (see the dots
surrounded by small circles at the bottom right of the picture).}
\label{henon_fig2}
\end{figure}
\begin{figure}[ht]
\begin{center}
\includegraphics[width=13cm,height=8cm]{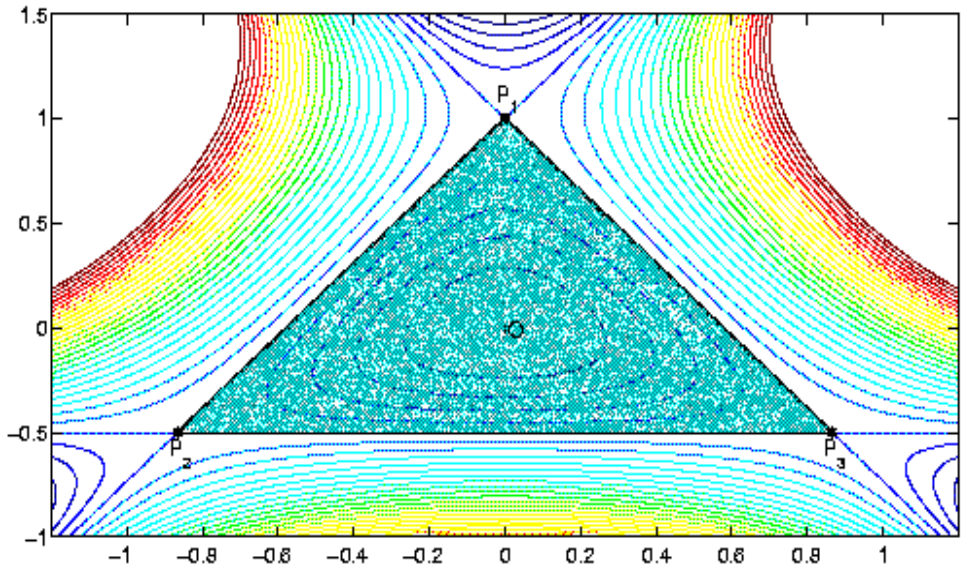}
\end{center}
\vspace*{-0.7cm}  \caption{The numerical trajectory in the
$(q_1,q_2)$-plane computed by the HBVM($4$,$2$) method  with
stepsize $h=1$. Since this method precisely conserves the total
energy of the system, the orbit is entirely contained in the
triangle at all times.} \label{henon_fig3}
\end{figure}

We have integrated problem \eqref{HH} in the time interval
$[0,5\cdot 10^4]$ with stepsize $h=1$ by using the Gauss method of
order four (HBVM($2$,$2$)) and the HBVM($4$,$2$) method which
assures an exact conservation of the total energy.

\begin{figure}[ht]
\begin{center}
\includegraphics[width=13cm,height=8cm]{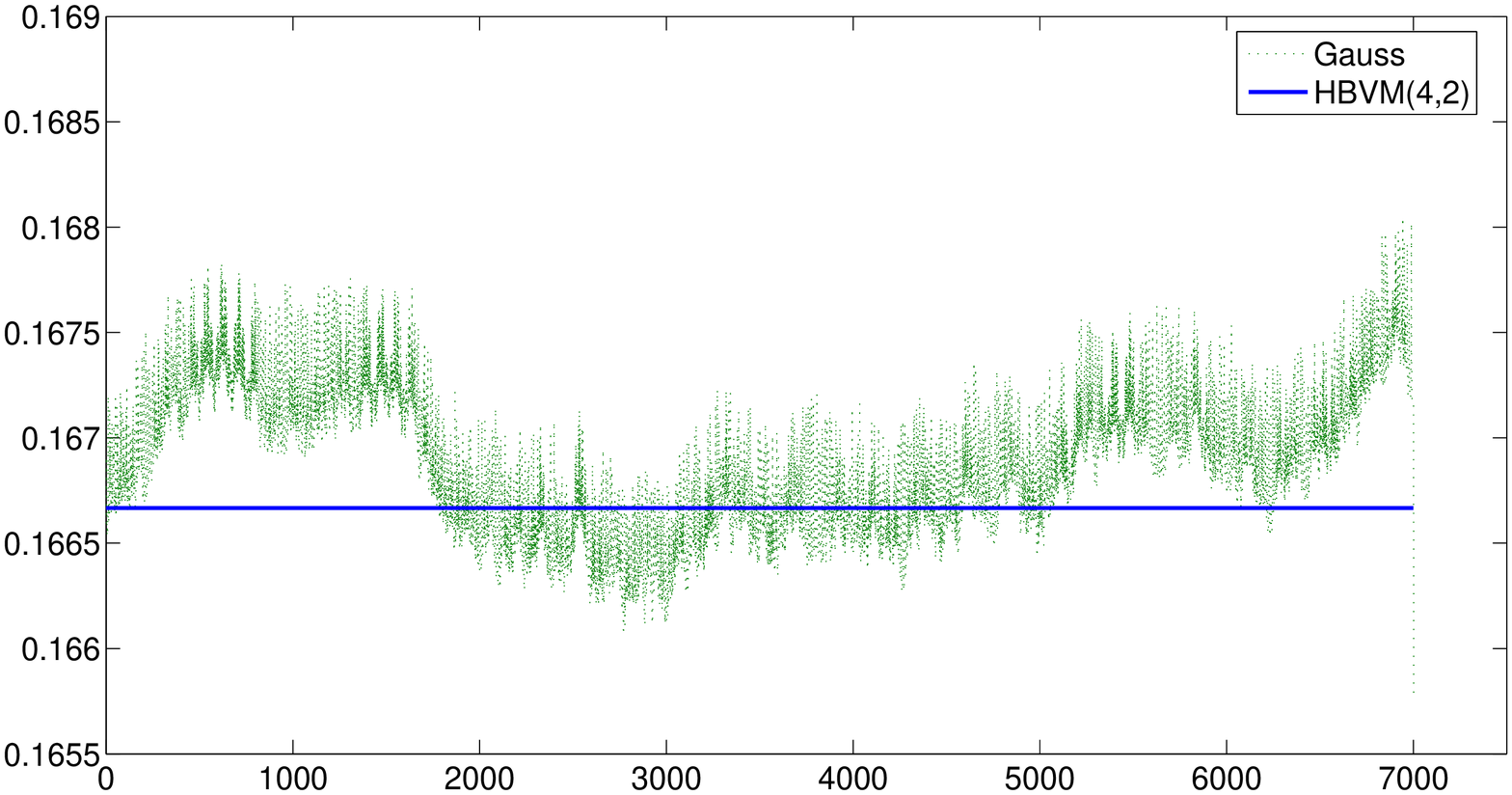}
\end{center}
\vspace*{-0.7cm}  \caption{Hamiltonian function evaluated along the
numerical solution of the Gauss and HBVM($4$,$2$) methods. The
irregular oscillations introduced by the Gauss method will cause the
associated numerical solution to eventually leave the stability
region centered at the origin.} \label{henon_fig4}
\end{figure}

Figures \ref{henon_fig2} and \ref{henon_fig3} show the numerical
trajectories in the $(q_1,q_2)$-plane as dots that eventually will
densely fill the triangle. The orbit generated by the Gauss method
is plotted up to time $t\simeq 7000$, since it then escapes from the
triangle, as highlighted by the three circles close to the saddle
point $P_3$. In fact, as Figure \ref{henon_fig4} shows, the
numerical Hamiltonian function associated with the Gauss method
produces very irregular oscillations around the theoretical value
(straight line) which eventually determine a loss of stability.

On the contrary, all the $50000$ dots of the numerical trajectory
computed by the HBVM($4$,$2$) method are visible  in Figure
\ref{henon_fig3}.

\subsection{Computing the period annulus of a non-degenerate center of
a polynomial Hamiltonian planar system.}

Non-degenerate centers\footnote{\,We recall that a center is an
equilibrium point which is surrounded by periodic orbits. It is
\textit{non-degenerate} if the linearized vector field at this
point has non-zero eigenvalues.} of planar, in particular
polynomial, Hamiltonian systems are extensively researched in the
modern literature (see \cite{DuLlAr,ChDe,LlRo,CiGaMa} and
references therein). The integration of such systems by means of
HBVMs deserves a particular interest because, the degrees of
freedom being one, the corresponding numerical solution is
guaranteed to lie on the same level set $H(q,p)=H(q_0,p_0)$ as the
theoretical orbit. Furthermore, if this latter consists of a
closed orbit surrounding an equilibrium point (center), the
numerical solution will (in general) fill densely the
corresponding closed level curve, thus reproducing the very same
phase portrait associated with the original continuous problem.

The region of marginal stability of a center $P_0$, is called the
\textit{period annulus of $P_0$} and will be denoted by $\cal P$:
it is the largest punctured neighborhood of the center consisting
of only periodic orbits.  The function which associates to any
periodic orbit in $\cal P$ its period is called the {\em period
function} of the center. Such function has been being intensively
studied for many years: its behavior relates to problems of
isochronicity,\footnote{\,Namely, all the orbits surrounding the
center $P_0$ share the same period.} monotonicity, bifurcation of
its critical points, etc.

The aim of the present example is to consider one such system and
try to reproduce numerically, as best as possible, the set $\partial
\cal P$, that is the boundary of the period annulus $\cal P$. Let
$H^\ast<+\infty$ be the value of the Hamiltonian function
corresponding to any points on $\partial \cal P$.\footnote{\,Here we
assume that the center $P_0$ is non global: this is certainly true
if $H(q,p)$ is a polynomial of odd degree.} The Hamiltonian function we consider
here is the fifth-degree
polynomial
\begin{equation}
\label{xavier} H(p,q)=A(p)+B(p)q+C(p)q^{2}+D(p)q^{3},
\end{equation}
where
$$
\begin{array}{lcl}
A(p)=p^{2}(\frac{1}{2}+c_{3}p+b_{3}p^{2}+a_{3}p^{3}), &\qquad&
B(p)=p^{2}(c_{2}+b_{2}p+a_{2}p^{2}), \\
C(p)=\frac{1}{2}+c_{1}p+b_{1}p^{2}+a_{1}p^{3}, &&
D(p)=c_{0}+b_{0}p+a_{0}p^{2},
\end{array}
$$
with $(a_0,a_1,a_2,a_3)\not = (0,0,0,0)$.\footnote{\,Otherwise the
degree of $H(q,p)$ becomes lower than $5$.} Note that, since
$H(q,p)=\frac{1}{2}(p^2+q^2) + \mbox{h.o.t.}$, we can assume $P_0$
to be the origin $O=(0,0)$.

The class of Hamiltonian systems defined by \eqref{xavier} has
been proposed in \cite{JaVi1} and \cite{JaVi2}.\footnote{\,The
authors showed that, without loss of generality, the form
\eqref{xavier} may be associated to any polynomial Hamiltonian
system of degree four and admitting a non-degenerate center, via a
suitable change of coordinates.} Their main result was  proving
that the origin may not be an isochronous center \cite{JaVi1} and,
more specifically, that the period tends to infinity as
$H(q_0,p_0) \nearrow H^\ast$, $(q_0,p_0)$ being the initial
condition associated with the differential system.

\begin{figure}[hb]
\begin{center}
\includegraphics[width=13cm,height=8cm]{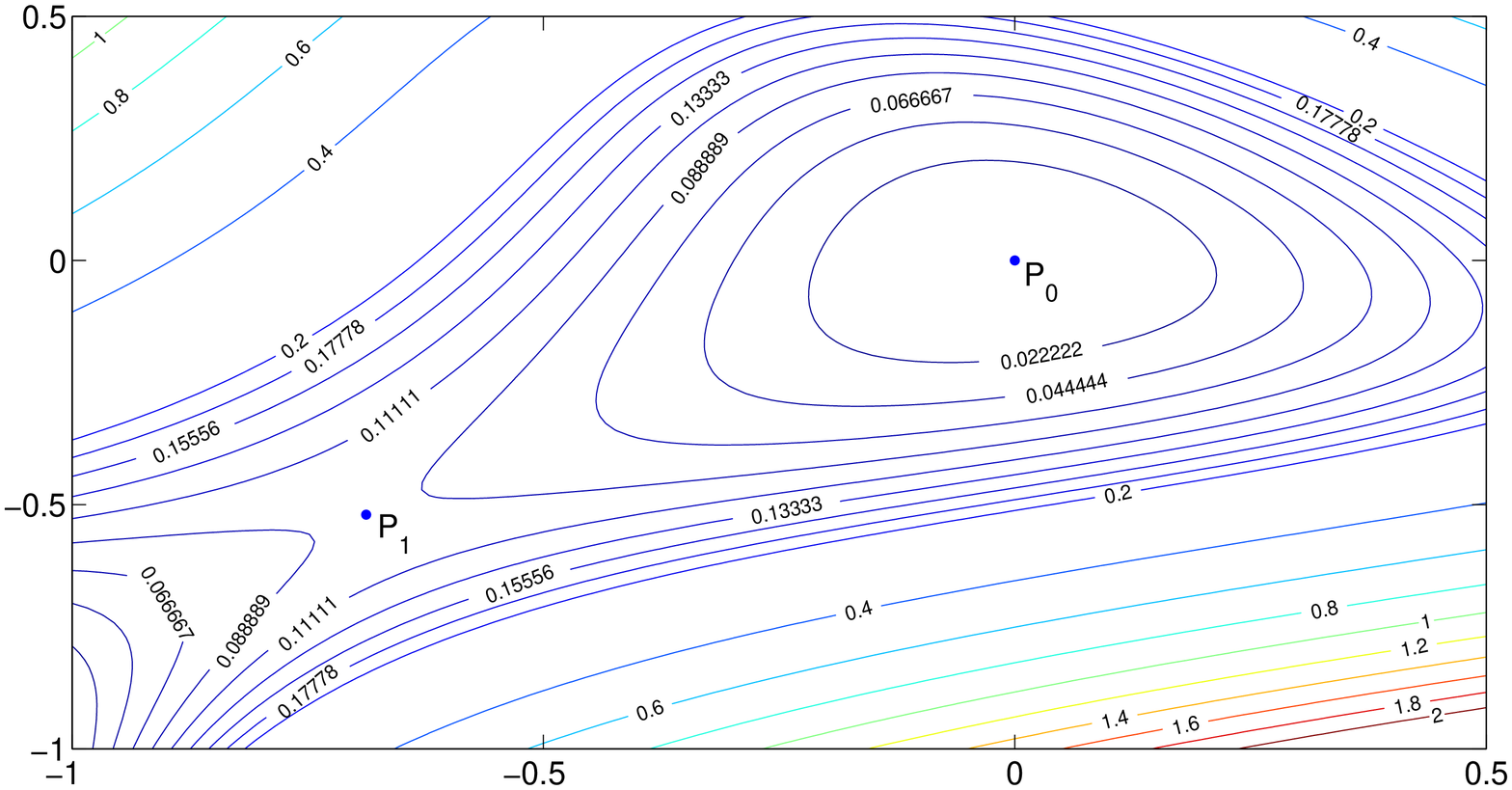}
\end{center}
\vspace*{-0.7cm}  \caption{Level curves of the Hamiltonian
\eqref{xavier} in a region that embraces the center point $P_0$
and the saddle point $P_1$. Each level curve, corresponding to an
orbit of the associated Hamiltonian system, is labeled by a number
that indicates its elevation.} \label{poinc_fig1}
\end{figure}

For our experiments, we have set the values of the coefficients
$\{a_i\}$, $\{b_i\}$, and $\{c_i\}$ as follows:
\begin{align}
a_0=0;  && a_1=0; && a_2=1; && a_3=0; \notag \\
b_0=0;  && b_1=1; && b_2=0; && b_3=1;  \label{poinc_data} \\
c_0=0;  && c_1=1; && c_2=1; && c_3=0 \notag.
\end{align}
In such a case, besides the origin $P_0=(0,0)$, $H(q,p)$ admits
the following real equilibrium points (up to the machine
precision):
$$\begin{array}{l}
P_1=(-6.879526475540134 \cdot 10^{-1},~ -5.206527058470621\cdot
10^{-1}) \longrightarrow \mbox{saddle point;} \\
P_2=(-1.179582379893681,~ 1.756351969248087) \longrightarrow \mbox{saddle point.} \\
\end{array}
$$
Figure \ref{poinc_fig1} reports the shape of the level curves of
\eqref{xavier}--\eqref{poinc_data} in a region enclosing $P_0$ and
$P_1$. We see that the limit closed orbit corresponding to
$\partial \cal P$ is the one embracing $P_0$ and having $P_1$ as
both $\omega$-limit point and $\alpha$-limit point\footnote{\,That
is, $\lim_{t\rightarrow \pm \infty} (q(t),p(t)) = P_1$ for any
choice of $(q_0,p_0) \in \partial P$.} and, therefore, the value
$H^\ast$ may be computed with precision as
\begin{equation}
\label{Hstar}
H^\ast = H(P_1)= 9.050199350868576 \cdot 10^{-2}.
\end{equation}

Now suppose we do not know the value $H^\ast$ in \eqref{Hstar} (it
will be used as a reference value)  and that we want to reproduce
the orbit covering $\partial \cal P$ by simply picking initial
points $(q_0,p_0)$ further and further away from the origin, and
checking whether the numerical solution remains bounded over  a
long time.\footnote{\,Of course, we cannot assume $(q_0,p_0)=P_1$
since $P_1$ is an equilibrium point.} More precisely, we will
locate the limit cycle by means of a dichotomic search, according
to the following algorithm:
\begin{itemize}
\setlength{\itemsep}{0cm} \setlength{\parskip}{0cm}
\item[step 1:] find a point $Q$ from which an orbit originates that
does not embraces the critical point $P_0$ (that is $Q \not \in \cal
P$);
\item[step 2:] consider the segment joining $P_0$ to $Q$:\\
\centerline{$\gamma(c)=(1-c)P_0+cP_1, \qquad c\in [0,1],$} and set
$c_0=0$ and $c_1=1$;
\item[step 3:] if $c_1-c_0<\mathrm{tol}$, STOP \qquad ($\mathrm{tol}$
is a specified tolerance);
\item[step 4:] set $c=\frac{c_0+c_1}{2}$ and solve numerically the
Hamiltonian problem defined in
\eqref{xavier}, considering $\gamma(c)$ as initial condition, in the
time interval $[0, hN]$ where $h>0$ is the stepsize and $N$ is a
positive integer such that $hN$ is large enough to give some
information about the fate of the orbit originating from
$\gamma(c)$.
\item[step 5:] if the numerical solution eventually depart from
$P_0$, set $c_1=c$, otherwise set $c_0=c$, go to step 3;
\end{itemize}

The point $y_0 \equiv (q_0,p_0)= \gamma(c)$, where $c$ is the
value resulting after the execution of the above procedure,  may
be assumed as a point on $\partial \cal P$ within the specified
tolerance $\mathrm{tol}$. Detecting the limit cycle with high
accuracy requires a huge number of simulations and therefore large
run times, also taking into account the wide time intervals that
must be used in order to inspect the asymptotic behavior of the
numerical solution.\footnote{Actually, by virtue of their
conservation properties, HBVMs do not need to be integrated over a
long time, even though here we do that for comparison purposes.}
Consequently, it would be advisable to work with a relatively
large stepsize $h$. We have set:
$$
h=1,~0.5, \qquad N=2500,~5000, \qquad \mathrm{tol}=2^{-52}~ \mbox{(i.e., the
value of {\tt eps} in Matlab)}, \qquad  Q=(0,1),
$$
to cover the integration interval $[0,2500]$.

\begin{figure}[ht]
\begin{center}
\includegraphics[width=7.7cm,height=6.5cm]{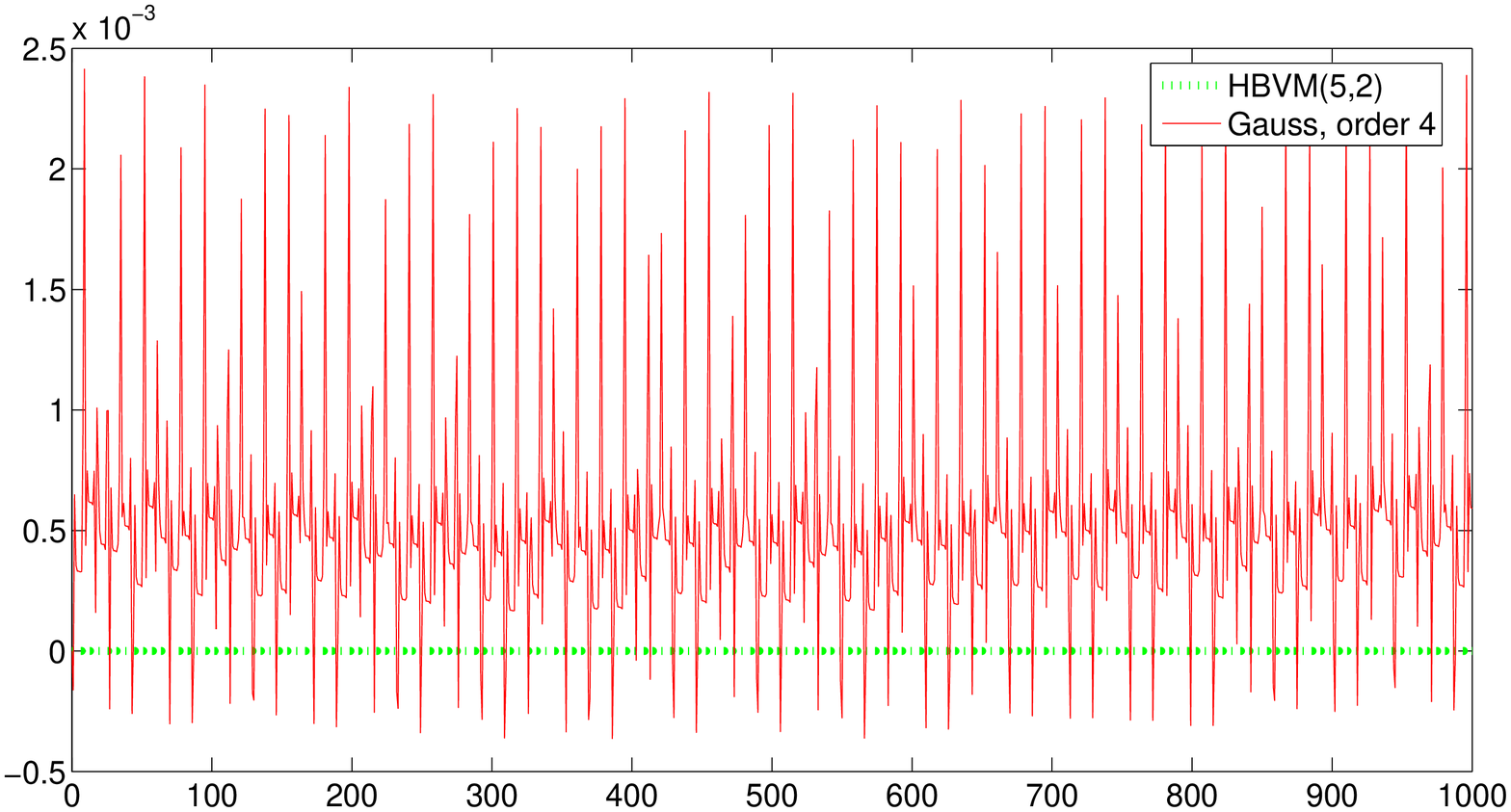}
\includegraphics[width=7.7cm,height=6.5cm]{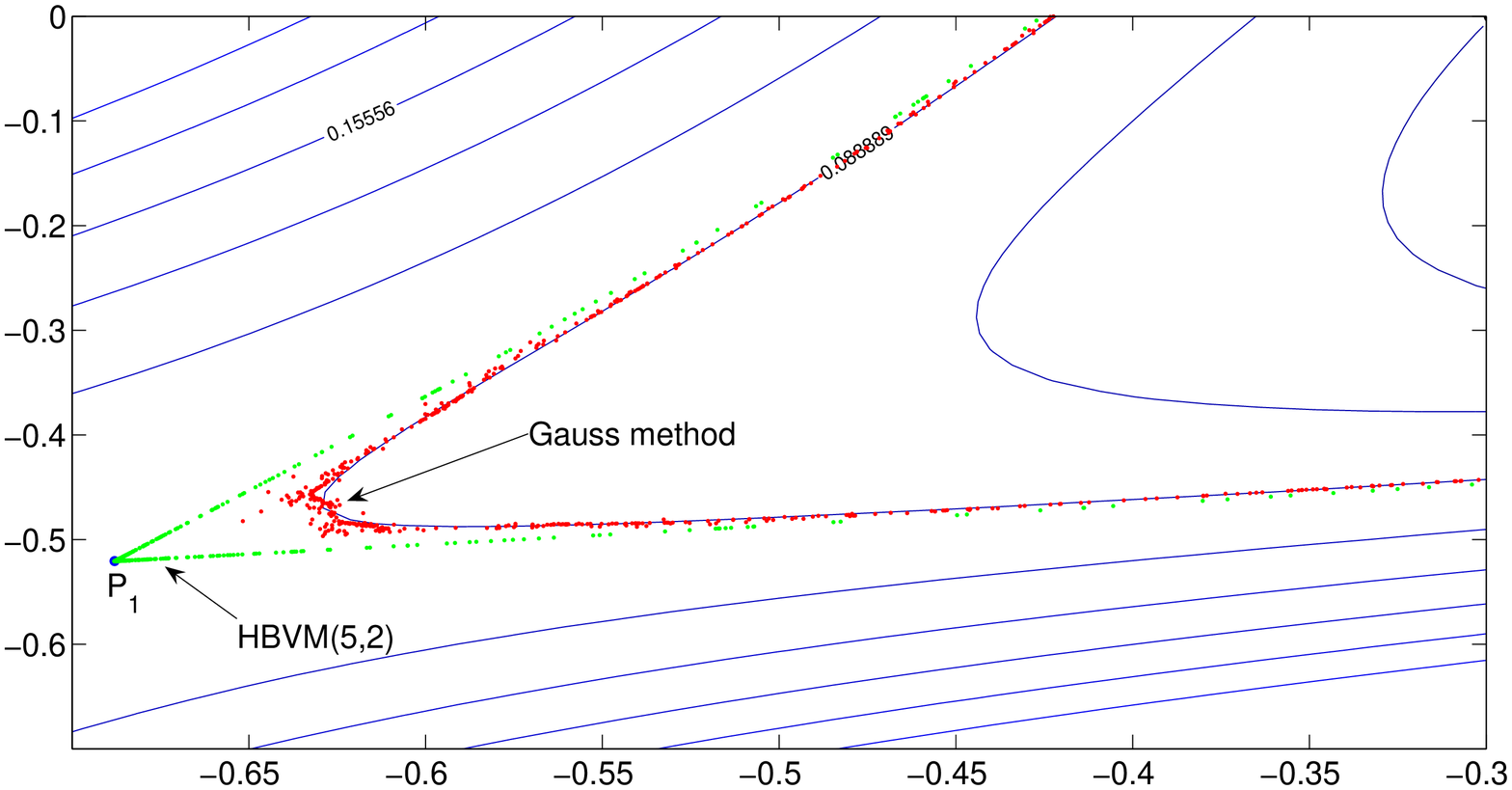}
\end{center}
\vspace*{-0.7cm} \caption{Left picture: Error
$H(q_n,p_n)-H(q_0,p_0)$ in the Hamiltonian function corresponding
to the numerical solutions computed by the Gauss method of order
$4$ and HBVM(5,2) (order 4), with stepsize $h=1$ and initial
conditions $y_0^{(2,2)}$ and $y_0^{(5,2)}$ respectively. Right
picture: a closeup of the two numerical orbits in a neighborhood
of the saddle point $P_1$ reveals the difficulty of  the Gauss
method in detecting the boundary of the period annulus.}
\label{poinc_fig2}
\end{figure}

\begin{table}[ht]
{\small $$
\begin{array}{|c|c|cccc|}
\hline
h &s & \begin{array}{c} \mbox{a point } y_0^{(s,s)}\in \partial {\cal P} \mbox{ computed} \\[-.12cm] \mbox{by the Gauss method} \end{array} & \frac{|H(y_0^{(s,s)})-H^\ast|}{H^\ast} & \begin{array}{c} \mbox{a point } y_0^{(k,s)} \in \partial {\cal P} \mbox{ computed} \\[-.12cm] \mbox{by HBVM($k$,$s$)} \end{array} & \frac{|H(y_0^{(k,s)})-H^\ast|}{H^\ast} \\
\hline
&2&(0,\,  3.723580509957994\cdot 10^{-1}) & 2.15\cdot 10^{-2}   &         (0,\,  3.757055929263451\cdot 10^{-1}) & 7.66\cdot 10^{-16} \\
\displaystyle  1 &3&(0,\,  3.748759009745006\cdot 10^{-1}) & 5.38\cdot 10^{-3}   &         (0,\,  3.757055929263451\cdot 10^{-1}) & 4.60\cdot 10^{-16} \\
&4&(0,\,  3.754691919292651\cdot 10^{-1}) & 1.53\cdot 10^{-3}   &         (0,\,  3.757055929263450\cdot 10^{-1}) & 1.22\cdot 10^{-15} \\
&5&(0,\,  3.756914213384024\cdot 10^{-1}) & 9.20\cdot 10^{-5}   &         (0,\,  3.757055929263451\cdot 10^{-1}) & 4.60\cdot 10^{-16} \\
\hline
&2&(0,\,  3.756045691696934\cdot 10^{-1}) & 6.56\cdot 10^{-4}   &         (0,\,  3.757055929263451\cdot 10^{-1}) & 4.60\cdot 10^{-16} \\
\displaystyle \frac{1}{2}&3&(0,\,  3.756828289241957\cdot 10^{-1}) & 1.47\cdot 10^{-4}   &         (0,\,  3.757055929263451\cdot 10^{-1}) & 4.60\cdot 10^{-16} \\
&4&(0,\,  3.757049796804918\cdot 10^{-1}) & 3.98\cdot 10^{-6}   &         (0,\,  3.757055929263451\cdot 10^{-1}) & 4.60\cdot 10^{-16} \\
&5&(0,\,  3.757055571549585\cdot 10^{-1}) & 2.32\cdot 10^{-7}   &         (0,\,  3.757055929263451\cdot 10^{-1}) & 4.60\cdot 10^{-16} \\
\hline
\end{array}
$$}
\caption{A point $y_0$ on the boundary of the period annulus $\cal
P$ is computed by the Gauss and HBVM methods of orders $4$, $6$,
$8$ and $10$ (corresponding to $s=2,\,3,\,4,\,5$ respectively). By
their very nature, if used with a sufficient number of silent
stages, HBVMs produce a numerical orbit that precisely lie on the
same level set $H(q,p)=H(q_0,p_0)$ as the theoretical one,
therefore we see that HBVMs can locate the point $y_0$ with
extreme precision, whatever the order and/or the stepsize used. On
the contrary, Gauss methods produce a certain error that may be
lowered by reducing the stepsize of integration $h$ and/or by
raising their order.} \label{poinc_tab1}
\end{table}
Table \ref{poinc_tab1} compares the results obtained by using the
Gauss (HBVM($s$,$s$)) and HBVM($k$,$s$) methods of orders
$4,\,6,\,8$ and $10$ (therefore, since $s=2,3,4,5$, we must
choose, according to (\ref{kk}), $k=5,8,10,13,$ respectively, in
order for the HBVM($k$,$s$) to exactly conserve the Hamiltonian
function). We have denoted by $y_0^{(k,s)}$ the point computed by
the method HBVM($k$,$s$), and reported the error
$|H(y_0^{(k,s)})-H^\ast|/H^\ast$ to estimate the accuracy with
which each method computes the boundary of $\cal P$. As was
expected,  the accuracy in detecting the right boundary of the
period annulus by means of HBVMs is of the same order as the
machine precision whatever the order and stepsize used (indeed,
the value of $y_0^{(k,s)}$ remains the same for all simulations).
On the contrary, the Gauss methods produce a certain error which
depends both on the stepsize and on the order used: increasing the
accuracy would require a suitable reduction of the stepsize and/or
a grow-up of the order. Figure \ref{poinc_fig2} shows that even
small oscillations of the numerical Hamiltonian function (left
picture) could produce a noticeable irregularity of the numerical
orbit in a neighborhood of the boundary of the period annulus
(right picture). By their very nature, HBVMs succeed in detecting
the set $\partial \cal P$ with an accuracy of the same order as
the machine precision: the error in the Hamiltonian function is
negligible (left picture) and the numerical  orbit correctly
passes through the saddle point $P_1$.

\end{document}